\documentclass[11pt]{amsart}
\usepackage{amsfonts,amsmath,amsthm,amscd,amssymb,latexsym,cite,verbatim,texdraw,floatflt,
caption2}
\RequirePackage{hyperref}
\usepackage[a4paper]{geometry}

\title{Selectors and orderings of  coarse spaces   }
\author{ Igor  Protasov}

\address{I.Protasov: Taras Shevchenko National University of Kyiv, Department of Computer Science and Cybernetics, Academic Glushkov pr. 4d, 03680 Kyiv, Ukraine}
\email{i.v.protasov@gmail.com}

\begin{document}
\begin{abstract} 
Given a coarse space 
$(X, \mathcal{E})$, we consider linear orders on 
$X$ compatible with the coarse structure  $\mathcal E$ 
and explore interplays between these orders and macro-uniform selectors of $(X, \mathcal{E})$.

\end{abstract}
\maketitle



\section{ Introduction and preliminaries}

The notion of selectors comes from {\it Topology}.
Let $X$ be a topological space, $exp \ X$ denotes the set of all non-empty closed subsets of $X$ endowed with some (initially, the Vietoris) topology, $\mathcal{F}$ be a  non-empty closed subset of $exp \ X$. A continuous mapping $f: \mathcal{F}\rightarrow X$ is called an $\mathcal{F}$-selector of $X$ if $f(A)\in A$ for each $A\in \mathcal{F}.$

Formally, coarse spaces,  introduced independently in 
\cite{b9} and \cite{b13} can be considered as asymptotic counterparts of uniform topological spaces.
But actually, this notion is rooted in {\it Geometry, Geometrical Group Theory} and {\it Combinatorics},
see \cite{b13}, \cite{b3}, \cite{b5} and \cite{b9}.

The investigation of selectos of coarse spaces was initiated in  \cite{b8}. We begin with some basic definitions. 

\vspace{5 mm}

Given a set $X$, a family $\mathcal{E}$  of subsets of $X\times X$ is called a
{\it  coarse structure} on $X$ if

\begin{itemize}
\item{} each $E \in \mathcal{E}$  contains the diagonal $\bigtriangleup _{X}:=\{(x,x): x\in X\}$ of $X$;
\vspace{3 mm}

\item{}  if  $E$, $E^{\prime} \in \mathcal{E}$  then  $E \circ E^{\prime} \in \mathcal{E}$  and
$ E^{-1} \in \mathcal{E}$,    where  $E \circ E^{\prime} = \{  (x,y): \exists z\;\; ((x,z) \in E,  \ (z, y)\in E^{\prime})\}$,    $ E^{-1} = \{ (y,x):  (x,y) \in E \}$;
\vspace{3 mm}

\item{} if $E \in \mathcal{E}$ and  $\bigtriangleup_{X}\subseteq E^{\prime}\subseteq E$  then  $E^{\prime} \in \mathcal{E}$.
\end{itemize}

Elements $E\in\mathcal E$ of the coarse structure are called {\em entourages} on $X$.

For $x\in X$  and $E\in \mathcal{E}$ the set $E[x]:= \{ y \in X: (x,y)\in\mathcal{E}\}$ is called the {\it ball of radius  $E$  centered at $x$}.
Since $E=\bigcup_{x\in X}( \{x\}\times E[x]) $, the entourage $E$ is uniquely determined by  the family of balls $\{ E[x]: x\in X\}$.
A subfamily ${\mathcal E} ^\prime \subseteq\mathcal E$ is called a {\em base} of the coarse structure $\mathcal E$ if each set $E\in\mathcal E$ is contained in some $E^\prime \in\mathcal E^\prime$.

The pair $(X, \mathcal{E})$  is called a {\it coarse space}  \cite{b13} or  a {\em ballean} \cite{b9}, \cite{b12}. 


A coarse  spaces
$(X, \mathcal{E})$ is  called 
{\it connected} if, for any $x, y \in X$, there exists $E\in \mathcal{E}$ such that $y\in E[x]$. 
A subset  $Y\subseteq  X$  is called {\it bounded} if $Y\subseteq E[x]$ for some $E\in \mathcal{E}$,
  and $x\in X$.
If  $(X, \mathcal{E})$ 
is connected  then 
the family $\mathcal{B}_{X}$ of all bounded subsets of $X$  is a bornology on $X$.
We recall that a family $\mathcal{B}$  of subsets of a set $X$ is a {\it bornology}
if $\mathcal{B}$ contains the family $[X] ^{<\omega} $  of all finite subsets of $X$
 and $\mathcal{B}$  is closed   under finite unions and taking subsets. A bornology $\mathcal B$ on a set $X$ is called {\em unbounded} if $X\notin\mathcal B$.
A subfamily  $\mathcal B^{\prime}$ of $\mathcal B$ is called a base for $\mathcal B$ if, for each $B \in \mathcal B$, there exists $B^{\prime} \in \mathcal B^{\prime}$ such that $B\subseteq B^{\prime}$.

Each subset $Y\subseteq X$ defines a {\it subspace}  $(Y, \mathcal{E}|_{Y})$  of $(X, \mathcal{E})$,
 where $\mathcal{E}|_{Y}= \{ E \cap (Y\times Y): E \in \mathcal{E}\}$.
A  subspace $(Y, \mathcal{E}|_{Y})$  is called  {\it large} if there exists $E\in \mathcal{E}$
 such that $X= E[Y]$, where $E[Y]=\bigcup _{y\in Y} E[y]$.

Let $(X, \mathcal{E})$, $(X^{\prime}, \mathcal{E}^{\prime})$
 be  coarse spaces. 
 A mapping $f: X \to X^{\prime}$ is called
  {\it  macro-uniform }  if for every $E\in \mathcal{E}$ there
  exists $E^{\prime}\in \mathcal{E}^{\prime}$  such that $f(E(x))\subseteq  E^{\prime}(f(x))$
    for each $x\in X$.
If $f$ is a bijection such that $f$  and $f ^{-1 }$ are macro-uniform, then   $f  $  is called an {\it asymorphism}.
If  $(X, \mathcal{E})$ and  $(X^{\prime}, \mathcal{E}^{\prime})$  contain large  asymorphic  subspaces, then they are called {\it coarsely equivalent.}

For a coarse space 
$(X,\mathcal{E})$, 
we denote by 
$exp \ X$
 the family  of all non-empty  subsets of $X$
  and by $exp \ \mathcal{E}$ the coarse structure on $exp \ X$
  with the base $\{ exp \ E : E\in  \mathcal{E}\}$, where
 $$(A,B)\in exp \ E 
  \Leftrightarrow A \subseteq E[B], \ \ B\subseteq E[A],$$
and say that $(exp \ X, exp \ \mathcal{E} )$ is the {\it hyperballean} of 
$(X,\mathcal{E})$. 
For hyperballeans, see \cite{b4}, \cite{b10}, \cite{b11}.

Let  $\mathcal{F}$ be 
  a non-empty  subspace  of $exp \ X$. 
 We say that a macro-uniform mapping 
 $f: \mathcal{F} \longrightarrow X$ 
 is an $\mathcal{F}$-{\it selector} 
 of $(X,\mathcal{E})$ if $f(A)\in A$ for each $A\in \mathcal{F}$. In the case $\mathcal{F}\in [X]^2$,
 $\mathcal{F}= \mathcal{B}_X$
 and $\mathcal{F}= exp \  X$, an $\mathcal{F}$- selector is called a $2$-{\it selector}, 
 a {\it bornologous selector}  
 and a {\it global selector}  respectively.

  We recall that a connected coarse space $(X,\mathcal{E})$ is {\it discrete} if, for each $E\in  \mathcal{E}$, there exists a bounded subset $B$ of $(X,\mathcal{E})$ such that $E[x]=\{x\}$  for each $x\in X\setminus B$. Every bornology $\mathcal{B}$ on a set $X$ defines the discrete coarse space $X_\mathcal{B} = (X,\mathcal{E} _\mathcal{B})$, where  $\mathcal{E}_\mathcal{B}$ is a coarse structure with the base 
$\{ E_B: B\in  \mathcal{B}\}$, $E_B [x]=B$ if $x\in B$ and 
$E_B [x]= \{x\}$
if $x\in X\setminus B$. On the other hand,  every discrete coarse space $(X, \mathcal{E})$ coincides with 
$X_\mathcal{B}$, where $\mathcal{B}$ is the bornology of bounded subsets of $(X, \mathcal{E})$.

\vspace{7 mm}

{\bf Theorem 1 [8]. } 
{\it For a bornology $\mathcal{B}$ on a set $X$,  the discrete coarse space 
 $X_{\mathcal{B}}$
  admits  a 2-selector if and only if
   there exists a linear order $\leq$ on $X$
 such that the family of intervals
  $\lbrace [a,b]: a,b \in  X, a\leq b
  \rbrace$ is a  base for $\mathcal{B}$. 
 \vspace{7 mm}
}

 In section 2, we analyze  interrelations between linear orders compatible with coarse structures and selectors. In Section 3, we apply obtained  results to characterize cellular ordinal 
coarse spaces which admit global  selectors. We conclude  with Section 4 on selectors of universal spaces.

\section{ Selectors and orderings }

{\bf Proposition 1.}
Let $(X, \mathcal{E})$
be a coarse space, 
$f: [X]^2 \rightarrow X$, $f(A)\in A$ for each 
$A\in [X]^2 $. Then the following statements are 
equivalent 

\vspace{5 mm}
{\it (i)} 
{\it $f$ is a 2-selector};

\vspace{3 mm}

{\it (ii)} {\it for every $E\in \mathcal{E}$, there exists $F\in \mathcal{E}$ such that $E\subseteq F$ and 
if
$\lbrace x,y\rbrace \in [X]^2$,  $f(\lbrace x,y\rbrace)=x$ 
$(f(\{x,y\})=y)$
and $y\in X\setminus F[x]$ then $f(\lbrace x^\prime,y\rbrace)=x^\prime$ 
$(f(\{x^\prime,y\})=y)$
 for each $x^\prime\in E[x]$.

\vspace{5 mm}

Proof.}  $(i) \Rightarrow(ii)$. 
Let $E= \mathcal{E}$.
Since $f$ is macro-uniform, there exists $F\in \mathcal{E}$, $F=F^{-1}$, $E\subseteq F$
such that, for any $(A, A^\prime) \in exp \ E$, 
we have  
$(f(A), f(A^\prime)) \in  F$. 
Let $A=\lbrace x, y \rbrace$, $A^\prime=\lbrace x^\prime, y \rbrace$, 
$x^\prime\in E[x]$, $f(\lbrace x, y\rbrace)=x$.
Then 
$f(\lbrace x, y\rbrace), f(\lbrace x^\prime, y\rbrace) \in F$, 
$(x,f(\lbrace x^\prime, y\rbrace) \in F$ so 
$f(\lbrace x^\prime, y\rbrace) = x^\prime$.
The case  $(f(\{x,y\})=y)$ is analogical.

\vspace{5 mm}

$(ii) \Rightarrow(i)$. Let $E\in \mathcal{E}$, 
$E= E^{-1}$ and let $F\in \mathcal{E}$, 
$F= F^{-1}$ is given by $(ii)$. To verify that $f$ is macro-uniform,  we show that if $A, A^\prime\in [X]^2$ 
and           $(A, A^\prime)\in exp \ E $
then             
$(f(A), f(A^\prime))\in F $.

Let $A=\lbrace x,y\rbrace$, 
$f(\lbrace x,y\rbrace)=x$,
$A^\prime=\lbrace x^\prime,y^\prime\rbrace, f(\lbrace x^\prime,y^\prime\rbrace) =x^\prime$.
We suppose that $(x,x^\prime)\notin F$ and 
$f(\lbrace x^\prime,x\rbrace) =x$.
By the choice of $F$, 
$f(\lbrace x^\prime,z\rbrace) =z$ for each 
$z\in E[x]$.
Since $E[x]\cap A^\prime \neq \emptyset$,
we have $y^\prime \in E[x]$ so 
$f(\{x^\prime, y^\prime \})= y^\prime$, 
contradicting 
$f(\{x^\prime, y^\prime \})= x^\prime$.
Hence, $(x,x^\prime)\in F$.
The case  $(f(\{x^\prime,x\})=x^\prime)$ is analogical.
 $ \ \  \  \Box $
 \vspace{7 mm}
 
Let $(X, \mathcal{E})$ be a coarse space.
We say that a linear order $\leq$ on $X$ is 
{\it compatible with the coarse structure} $\mathcal{E}$
if, for every $E\in \mathcal{E}$, there exists $F\in \mathcal{E}$ such that  
$E\subseteq F$ and if $\{x,y\}\in [X]^2$, $x< y$ 
$(y<x)$
 and $y\in X\setminus F[x]$  then $x^\prime <y$ 
$(y<x^\prime)$ 
  for each $x^\prime \in E[x]$.

\vspace{7 mm}
 
 {\bf Proposition  2. } 
{\it Let $(X, \mathcal{E})$ be a coarse space and let $\leq$ be a linear order on $X$ compatible with $\mathcal{E}$. Then the following statements hold

\vspace{3 mm}

$(i)$  the mapping  $f: [X]^2 \rightarrow X$, 
defined by $f(A)= min \ A$, is a 2-selector of $(X, \mathcal{E})$;

\vspace{3 mm}

$(ii)$ for every $E\in \mathcal{E}$, there exists $H\in \mathcal{E}$  such that $E\subseteq H$ and if $A, A^\prime \in [X]^2$ and $(A, A^\prime) \in exp \ E$ then 
$(min \ A, \ min \ A^\prime) \in H$; 
 
\vspace{3 mm}

$(iii)$  if $(X, \mathcal{E})$  is connected then, for any $a, b \in X$, $a<b$, the interval $[a,b]=\{ x\in X: a\leq x \leq b \}$ is bounded in $(X, \mathcal{E})$. 
 
 \vspace{5 mm}

 Proof. } 
The statement (i) follows from Proposition 1, $(ii)$ follows from $(i)$.

To prove $(iii)$, we use the connectedness of 
$(X, \mathcal{E})$ to find $E\in \mathcal{E}$, $E=E^{-1}$ such that $(a,b)\in E$.
Then we take $F\in \mathcal{E}$, $F= F^{-1}$
given by the definition of an order compatible with the 
coarse structure. We assume that $[a,b]$ is unbounded and choose $ c \in [a,b]$, 
$a< c <b$ such that $c\in X\setminus F[a]$.
Then $x<c$ for each $x\in E[a]$, in particular $b<c$ and we get a contradiction. $ \ \  \  \Box $
 \vspace{7 mm}

 {\bf Proposition 3. } 
{\it Let $(X, \mathcal{E})$ be a 
coarse space, $\leq$  be a well order on $X$
compatible with 
 $\mathcal{E}$. Then 
 $X_\mathcal{E}$ has a 
  global selector.
  
  \vspace{5 mm}

 Proof. } For each $A\in exp \ X$, we put 
 $f(A)=min \ A $ and note that $f$ is a global selector.
 $ \ \  \  \Box $
 \vspace{7 mm}

 {\bf Proposition 4. } 
{\it Let $(X, \mathcal{E})$ be a 
connected 
coarse space with the bornology $\mathcal{B}$
of bounded subsets, $X_\mathcal{B}$ denotes the discrete coarse space defined by $\mathcal{B}$.
If $f$ is a 2-selector of $(X, \mathcal{E})$ then 
$f$ is a  2-selector of $X_\mathcal{B}$.

 \vspace{5 mm}

 Proof. } For each $B\in \mathcal{B}$, we denote by 
 $E_B$  the set $\{(x.y): x,y\in B \}\cup \vartriangle_X$.
 Then 
 $\{E_B : B\in\mathcal{B}\}$ is the coarse  structure of $X_\mathcal{E}$ and 
 $E_B \in \mathcal{E}$ for each  $B\in\mathcal{B}$.
  
 Let $A, A^\prime \in [X]^2$ and 
 $(A, A^\prime) \ \in exp \ E_B$.
 Since $f$ is a 
  2-selector of $(X, \mathcal{E})$, there exists $F\in \mathcal{E}$, $F=F^{-1}$ such that $(f(A), f(A^\prime))\in F$.

If $A\cap B=\emptyset$ then $A= A^\prime$.
If $A\subseteq B$ then $A^\prime\subseteq B$, so $(f(A), f(A^\prime))\in E_B$. 

Let $A=\{b, a\}$,  $A^\prime=\{b^\prime, a\}$,
$b\in B$, 
$b^\prime\in B$ and  $a\in X\setminus B$. If $a\in F[\{b, b^\prime\}]$
 then $f(A), f(A^\prime)\in F[\{b, b^\prime\}] $.
 If $a\notin  F[\{b, b^\prime\}]$ then either 
 $f(A)= f(A^\prime) = a $ of 
 $f(A), f(A^\prime)\in \{b, b^\prime\} $.

In all considered cases, we have 
$(f(A), f(A^\prime))\in E_{F[B] }$.
Hence, $f$ is a   
2-selector of $X_\mathcal{B}$. $\ \  \  \Box $
 \vspace{7 mm}
 
 {\bf Proposition  5. } 
{\it Let $(X, \mathcal{E})$, $(X^\prime, \mathcal{E}^\prime)$ are coarsely equivalent.
If $(X^\prime, \mathcal{E}^\prime)$  admits a global 
selector then $(X, \mathcal{E})$ admits a global selector. 
The same is true for 2-selector and bornologous selectors. 
 \vspace{5 mm}

 Proof. }  We consider the case of global selector. Let 
 $f^\prime : exp \ X^\prime \rightarrow X^\prime$ is a 
 global selector of $(X^\prime), \mathcal{E}^\prime))$.
 We suppose that  
 $(X, \mathcal{E})$, $(X^\prime, \mathcal{E}^\prime)$ are asymorphic and 
 $h: (X, \mathcal{E})\rightarrow (X^\prime, \mathcal{E}^\prime)$ is an asymorphism. We denote by 
 $\overline{h}$ the natural extension $\overline{h}: exp \ X\rightarrow exp \ X^\prime$ of $h$.
 Then the straitforward verification gives that $h^{-1} f^\prime \overline{h}$ is a global selector of $(X, \mathcal{E})$.

\vspace{3 mm}

Now let $X^\prime$ is a large subset of  
$(X, \mathcal{E})$, $\mathcal{E}^\prime = \mathcal{E}_{X^\prime} $, $f^\prime : exp \ X^\prime\rightarrow X$ is a global selector of $(X^\prime, \mathcal{E}^\prime)$.
We take $H\in \mathcal{E}$ such that $X=H [X^\prime]$.
Let $Y\in exp \ X$. For each $y\in Y$, we pick $z_y\in X^\prime$  such that $y\in H[z_y]$. Let $Z=\{ z_y : y\in Y \}$
and $z=f^\prime (Z)$.
We take $x_z\in Y$ such that $x_z \in H[z]$ and put
$f(Y)=x_z$.
Then the straightforward verification gives us $f: exp \ X\rightarrow X$ is a global selector of $(X, \mathcal{E})$. $ \ \  \  \Box $
 \vspace{7 mm}

{\bf Question 1. } 
{\it Let $\leq$ be a linear order on $X$ 
compatible with $\mathcal{E}$. Is $\mathcal{E}$ an interval  coarse structure?}

\vspace{7 mm}

{\bf Question 2. } 
{\it Let a coarse space 
$(X, \mathcal{E})$ admits a global selector. Does there exist a linear order on $X$ compatible with $\mathcal{E}$?  }

\vspace{7 mm}

{\bf Question 3. } 
{\it Let a coarse space $(X, \mathcal{E})$ admits a 2-selector.
Does $(X, \mathcal{E})$  admit  a bornologous selector?}

\vspace{7 mm}

{\bf Question 4. } 
{\it Let a coarse space $(X, \mathcal{E})$ admits a 
  bornologous selector. 
Does $(X, \mathcal{E})$  admit a global selector? }

\section{ Selectors of  cellular spaces}

Let $(X, \mathcal{E})$
be a coarse space.
An entourage $E\in \mathcal{E}$ 
is called {\it cellular} if $E$ is an equivalence relation.
If $(X, \mathcal{E})$ is connected and $ \mathcal{E}$
 has a base consisting of cellular entourages then 
 $(X, \mathcal{E})$ is called cellular.
 By [12, Theorem 3.1.3], $(X, \mathcal{E})$ is cellular if and only if $asdim \ (X, \mathcal{E})=0$.
 
 Every discrete coarse space and every coarse space of an ultrametric space are cellular.
 
 Following [12, p.63], we say that a coarse space 
 $(X, \mathcal{E})$ is {\it ordinal} if $ \mathcal{E}$
 has a base well-ordered by inclusion. 
 We note that if $ \mathcal{E}$ has  a base linearly
 ordered by inclusion then $(X, \mathcal{E})$ is cellular. For the structure of cellular ordinal spaces, see 
 \cite{b1}.

 \vspace{5 mm}

 Let $\kappa$, $\gamma$ be cardinals. 
 Following  \cite{b1},
 we denote 
 $\kappa^{<\gamma}= \{ (x_\alpha)_{\alpha<\gamma} : 
 x_\alpha \in \kappa, \ 
 x_\alpha =0 \ $ 
 for all but finitely many $ \ \  \alpha<\gamma \}$,
$K_{\alpha}= \{ ((x_\alpha)_{\alpha<\gamma}, \ (y_\alpha)_{\alpha<\gamma}) : x_\beta =y_\beta$
 for each  $\beta\leq \alpha \}$.

We take the coarse structure $\mathcal{K}_\gamma$ with the base $\{K_\gamma: \alpha<\gamma \}$
and observe that each entourage $K_\alpha$ is cellular.
Thus,  the macrocube $(\kappa_\gamma , \mathcal{K}_\gamma)$ is cellular and ordinal. 

We denote ${\bf 0} = (x_\alpha)$,  $x_\alpha =0$
for each $\alpha<\gamma$ and, for $x=(x_\alpha)_{\alpha<\gamma}, \ x\neq 0$,
$max \ x= \{ max \ \alpha : x_\alpha \neq 0 \}$.
Given any 
$x=(x_\alpha)_{\alpha<\gamma}, \ 
y= (y_\alpha)_{\alpha<\gamma}, \ x\neq {\bf 0} \ y\neq 0$, we write $x\prec y$
if either $max \ x < max \ y$ or $max \ x = max \ y=\alpha$ and $x_\alpha < y_\alpha$.
Also, {\bf 0} $\prec x$ for $x\neq $ {\bf 0}.
Then $\preceq$ is a total order on $\kappa_\gamma$ 
compatible with the coarse structure $\mathcal{K}_\gamma$.

\vspace{7 mm}
{\bf Theorem 2. } {\it Every cellular ordinal space 
$(X, \mathcal{E})$
admits a  well-ordering compatible with $\mathcal{E}$.

\vspace{3 mm}

 Proof. } We put $\kappa=|X|$.
 By [1, Lemma 5.1], there exists an asymorphic embedding $f: (X, \mathcal{E})\rightarrow (\kappa, \mathcal{K}_\kappa)$.
 The total order $\preceq$ defined above on $\kappa^\kappa$ induces the 
total order $\preceq_{f(X)}$ on $f(X)$ 
compatible with the coarse of the subspace 
$f(X)$ of $(\kappa, \mathcal{K}_\kappa)$.
Applying $f^{-1}$,
we get the desired total order on $(X, \mathcal{E})$.
$ \ \  \  \Box $
 \vspace{7 mm}

{\bf Theorem 3. } {\it Every cellular ordinal space 
$(X, \mathcal{E})$
admits a  global selector
 $f: exp \ X \rightarrow X$.

\vspace{3 mm}

 Proof. } Apply Theorem 2 and Proposition 3. 
 $ \ \  \  \Box $
 \vspace{7 mm}

{\bf Question 5. } {\it How can one detect whether a given 
cellular
coarse space admits a global selector?}
 
  \vspace{7 mm}

Now we apply obtained results to coarse spaces of groups. Let $G$ be a group with the identity. 
We denote by $\mathcal{E}_G$ the coarse structure of 
$G$ with the base  
$$\{\{ (x, y)\in G\times G : y\in Fx: F\in [G]^{<\omega},  \  e\in F\}$$ and say that $(G, \mathcal{E}_G)$
 is the {\it finitary coarse space } of $G$.
 It should be mentioned that finitary coarse spaces of groups are used as tools in {\it Geometric Group Theory}, see 
 \cite{b3},  \cite{b5}.

 \vspace{7 mm}
 {\bf Theorem 4. } {\it If a group $G$ is uncountable then  
$(G, \mathcal{E}_G)$
does not admit a
2-selector.

\vspace{3 mm}

 Proof. } We note that the bornology of bounded subsets 
 of $(G, \mathcal{E}_G))$ is $[G]^{<\omega}$.
 Apply Proposition 4  and Theorem 1. $ \ \  \  \Box $
 \vspace{7 mm}
 
 It is easy to see that 
$(G, \mathcal{E}_G))$ is cellular if and only if $G$ is locally finite (i.e. each finite subset of $G$ generates a finite subgroup. 

\vspace{7 mm}

 {\bf Theorem 5. } {\it If $G$  is a  countable locally finite  group then  the  finitary coarse space 
$(G, \mathcal{E}_G)$ admits a global 
selector. \vspace{3 mm}

 Proof. } We note that $\mathcal{E}_G))$
 has a countable base and apply Theorem 2.$ \ \  \  \Box $
 \vspace{7 mm}
 
 Any two countable  locally finite groups are coarsely
 equivalent \cite{b2}, for classification of countable locally finite groups up to asymorphisms, see \cite{b6}.

\section{Selectors  of universal spaces}

 Let $X$ be a set, $\mathcal{E} \subseteq X\times X$,
  $\delta_X \subseteq X$.  We say that an entourage $E$ is 
  

\vspace{7 mm}

\begin{itemize}
\item{} {\it locally finite} if 
 $E [x]$,  $E^{-1} [x]$
 are finite for each  $x\in X$;
\vspace{3 mm}

\item{}  {\it  finitary} if there exists a natural  number $n$ such that 
$|E[x]|< n$, $|E^{-1}[x]| < n$  for each $x\in X$.

\end{itemize}

\vspace{5 mm}

A coarse space $(X, \mathcal{E})$ is called 
{\it locally finite (finitary)} if each entourage 
$E\in \mathcal{E}$ is locally finite (finitary). 
If $E, H$
are locally finite (finitary) then $E\circ H$, $E^{-1}$
are locally finite (finitary). 
We denote 
\vspace{4 mm}

$ \ \ \Lambda= \{ E: E\in \omega\times\omega, \ E $ is a locally finite entourage$\},$ 

$ \ \ \mathcal{F}= \{ E: E\in \omega\times\omega, \ E $ is a finitary entourage$\},$
\vspace{4 mm}

\noindent and say that $(\omega, \Lambda)$ $(resp. (\omega, \mathcal{F}))$
is the universal  {\it locally finite} (resp. {\it finitary}) space.

We denote by $S_\omega$ the group of all permutations of $\omega$, $id$ is the identity permutation.
By [7, Theorem 3], the coarse  structure $\mathcal{F}$ has the base $$\{\{ (x,y): x\in Fy\}: F \in [S_\omega]^{\omega}, \ \  id \in F \}. $$

{\bf Theorem 6.} {\it The coarse space $(\omega, \Lambda)$ admits a global selector.

\vspace{5 mm}

Proof.} We denote by $\leq$ the natural order on $\omega$, prove that $\leq$
is compatible with  $\Lambda$ and apply Proposition 3. 

For $E\in \Lambda$, let 
$\overline{E}= \{ (x, y): min E [x] \leq y \leq max E[x]\}. $
Clearly, $\overline{E}\in \Lambda$. If $x, y\in X$,
$x<y$ and $y\in \omega\setminus \overline{E}[x]$
then $x^\prime < y$ for each  $x^\prime \in E[x]$. $ \ \  \  \Box $
\vspace{7 mm}

{\bf Theorem 7.} {\it The coarse space $(\omega, \mathcal{F})$  does not admit 2-selectors.

\vspace{5 mm}

Proof.} We suppose the contrary and let $f$ be a 2-selector of $(\omega, \mathcal{F})$. We define a binary relation $ \prec$ on $X$ by $x\prec y$
if and only if  $x\neq y$ and $f(\{x,y\})=x$.
Then we choose inductively an injective sequence $(a_n)_{n\in\omega}$ in $\omega$ such that either $a_i\prec a_j$ for all  $i<j$ or $a_j \prec a_i$ for all $i<j$. 
We consider only the first case, the second is analogous. 

We partition $\{a_n : n<\omega \}$ into consequive with respect to $\prec$ intervals $\{ T_n : n<\omega\}$ of length $2n+1$.
We define a permutation $h$ of order 2 of $\omega$ as follows. 
For $x\in \omega \setminus \{ a_n: n<\omega \}$, $hx=x$.
We take $T_n$, $T_n = \{ a_m , \dots , a_{m+2n+1} \}$  
and put $h a_m =a_{m+2n+1}$, 
$h a_{m+1} =a_{m+2n}, \dots ,$
$h a_{m+n+1} =a_{m+n+1}$.
We put $F=\{h, id\}$,
$E=\{ (x,y): y\in Fx \}$.
Since $f$ is macro-uniform, there exists $H\in \mathcal{F}$ such that  if $A, A^\prime \in [X]^2$,
$A\subseteq E[A^\prime]$, $A^\prime\subseteq E[A]$
then $(f(A), f(A^\prime)) \in H$.

We take $k$ such that $|H[x]|< k$ for each $x\in \omega$. Let $n>k$. Since  
$(\{a_{m+i},  a_{m+n+1}\},$ 
$ \{a_{m+n+1}, a_{m+2n+1 -i })\in E$
for each $i\in \{0, \dots , n_1 \}$, 
$(a_{m+i},  a_{m+n+1})\in H$ contradicting 
$|H[a_{m+n+1}]|<k$.
$\  \  \Box $

\end{document}